\documentclass[letterpaper, 10 pt, conference]{ieeeconf}

\IEEEoverridecommandlockouts
\usepackage{amsmath}
\usepackage{amssymb}
\usepackage{graphicx}
\usepackage{epstopdf}
\usepackage{cite}
\usepackage{tikz}
\usetikzlibrary{arrows.meta,positioning,calc}

\newtheorem{theorem}{Theorem}
\newtheorem{lemma}{Lemma}
\newtheorem{assumption}{Assumption}
\newtheorem{definition}{Definition}
\newtheorem{remark}{Remark}

\newcommand{\sgn}{\operatorname{sgn}}
\newcommand{\R}{\mathbb{R}}
\newcommand{\N}{\mathbb{N}}

\title{\LARGE \bf
Small-Gain Analysis of Robust Event-Triggered\\ Extremum Seeking for Scalar Maps}

\author{Victor Hugo Pereira Rodrigues$^{1}$, Tiago Roux Oliveira$^{1}$, Miroslav Krsti\'c$^{2}$ and Zhong-Ping Jiang$^{3}$
\thanks{This work was financed in part by the Coordena\c{c}\~ao de Aperfei\c{c}oamento de Pessoal de N\'ivel Superior -- Brasil (CAPES) -- Finance Code 001, and by the Brazilian funding agencies CNPq and FAPERJ.}
\thanks{$^{1}$V. H. P. Rodrigues and T. R. Oliveira are with the Department of Electronics and Telecommunication Engineering, State University of Rio de Janeiro (UERJ), Rio de Janeiro, RJ, Brazil.
        {\tt\small victor.rodrigues@uerj.br, tiagoroux@uerj.br}}%
\thanks{$^{2}$M. Krsti\'c is with the Department of Mechanical and Aerospace Engineering, University of California San Diego (UCSD), La Jolla, CA, USA.
        {\tt\small mkrstic@ucsd.edu}}%
\thanks{$^{3}$Z.-P. Jiang is with the Department of Electrical and Computer Engineering, Tandon School of Engineering, New York University (NYU), Brooklyn, NY, USA.
        {\tt\small zjiang@nyu.edu}}%
}

\begin{document}

\maketitle
\thispagestyle{empty}
\pagestyle{empty}

\begin{abstract}
We study a sample-and-hold implementation of gradient extremum seeking for a scalar quadratic map in which the control signal is refreshed only when an event occurs. The event generator acts on the demodulated gradient signal as it is actually measured, rather than on its averaged version, and contains a strictly positive offset, following the small-gain philosophy for robust event-triggered control. Rather than invoking averaging theorems, we work directly with the exact dynamics of the parameter estimation error. In the scalar setting, the demodulated gradient has an explicit sinusoidal structure, which is exploited through a change of variables that absorbs the zero-mean oscillatory perturbation, together with a quadratic Lyapunov function and the persistence of excitation of the probing signal. For a sufficiently fast dither, the estimation error is shown to decay exponentially to a residual neighborhood of the origin whose size is quantified by the dither amplitude and frequency, and by the relative threshold and the offset of the event generator. All constants, including the input-to-state stability gain of the extremum seeking loop, the admissible dither frequency and the positive lower bound on the time between two consecutive events, are given in closed form, so that Zeno behavior is excluded. A simulation example illustrates the results.
\end{abstract}

\section{INTRODUCTION}

Extremum seeking (ES) is a non-model-based adaptive technique that drives the input of an uncertain map toward the point where the output is optimal, using nothing more than output measurements and a periodic probing signal \cite{AK:2003,KW:2000,K:2014}. After the stability proof in \cite{KW:2000}, the method has been extended in several directions, including multi-parameter and Newton-type schemes \cite{GKN:2012}, maps with delays and distributed-parameter actuation \cite{OKT:2017,TRO_book2022}, and applications such as source seeking with autonomous vehicles \cite{ZAGSK:2007}.

When the optimizer runs on an embedded or networked platform, refreshing the control input continuously or at a fixed fast rate squanders bandwidth and processor time. Event-triggered control replaces time-driven updates by state-driven ones: the actuator value is recomputed only when a prescribed triggering condition fails \cite{T:2007,HJT:2012}. The idea has been developed, among many other settings, for nonlinear systems with dynamic triggers \cite{ZLJ:2017}. Bringing event-triggered updates into ES is more recent: the gradient scheme of \cite{ROHDK:2025} was followed by a Newton-type version \cite{ROKT:2026}, by a discrete-time multivariable formulation \cite{ROKA:2026}, and by game-theoretic settings in which several players seek a Nash equilibrium \cite{ROKB:2026}. All of these works share one methodological choice: the trigger is designed on the \emph{average} model, and the behavior of the true oscillating loop is inferred afterwards from averaging theorems.

Here we take a different path, motivated by the small-gain approach to robust event-triggered control of nonlinear systems \cite{LJ:2015,LZJ:2020}, which rests on input-to-state stability (ISS) and on the small-gain theory for interconnections \cite{PJ:2006}. We regard the event-triggered ES loop as the feedback interconnection of two blocks: the dynamics of the parameter estimation error, and the event generator. The sampling error and the oscillations produced by the dither enter this interconnection as disturbances. We focus on the scalar (single-input) case, for which the demodulated gradient admits a fully explicit sinusoidal decomposition and all the constants of the analysis can be written in closed form. The contributions are:
\begin{itemize}
\item a static triggering rule with a positive offset, evaluated on the measured gradient estimate itself (not on its average), which tolerates the persistent oscillatory perturbations that are intrinsic to ES;
\item an ISS estimate of the ES dynamics with respect to the sampling error, derived on the exact error equation, without averaging, through the sign structure and the persistence of excitation of the demodulated gradient;
\item a small-gain condition coupling the ISS gain of the ES dynamics with the gain of the event generator, which yields local exponential practical stability and an explicit residual set;
\item a closed-form, uniform lower bound on the inter-event times, which rules out Zeno behavior.
\end{itemize}
The results are the scalar specialization of a multivariable theory developed by the authors; in the scalar case the extra assumptions on probing frequencies and gain structure disappear and every constant becomes transparent.

\textit{Notation:} $|\cdot|$ is the absolute value, $\sgn(\cdot)$ the sign function and $\N=\{0,1,2,\ldots\}$. For a function $f$ and an interval $I$, $\sup_{I}|f|:=\sup_{\tau\in I}|f(\tau)|$.

\section{PROBLEM STATEMENT}

We consider the scalar static map
\begin{align}
y(t)=Q(\theta(t))=Q^{\ast}+\frac{H^{\ast}}{2}\left(\theta(t)-\theta^{\ast}\right)^{2}\,, \label{eq:map}
\end{align}
where $\theta(t)\in\R$ is the input that can be manipulated in real time, $y(t)\in\R$ is the measured output, $\theta^{\ast}\in\R$ is the optimizer, $Q^{\ast}=Q(\theta^{\ast})$ is the optimal value and $H^{\ast}\in\R$ is the curvature (second derivative) of the map.

\begin{assumption}\label{ass:map}
The constants $\theta^{\ast}$, $Q^{\ast}$ and $H^{\ast}$ are unknown, but $H^{\ast}\neq 0$ and $\sgn(H^{\ast})$ is known.
\end{assumption}

Since any smooth map with a nondegenerate extremum agrees with (\ref{eq:map}) up to second order near $\theta^{\ast}$, the quadratic model captures the local behavior of general scalar ES problems while keeping the demodulated signal explicitly computable. When $H^{\ast}>0$ the task is to minimize $Q$, and when $H^{\ast}<0$ the task is to maximize it; in both cases $\theta(t)$ must be steered to a small vicinity of $\theta^{\ast}$. We write
\begin{align}
h:=|H^{\ast}|>0\,. \label{eq:h}
\end{align}

\section{EVENT-TRIGGERED EXTREMUM SEEKING}

\subsection{Probing, Demodulation and Error Variable}

\begin{figure}[!t]
\centering
\begin{tikzpicture}[>=Latex,font=\footnotesize,
  block/.style={draw,minimum height=6mm,minimum width=7mm,inner sep=2pt},
  sum/.style={draw,circle,minimum size=4mm,inner sep=0pt}]
\node[block] (K) at (0,0) {$k$};
\node[block,right=7mm of K] (int) {$\displaystyle\int$};
\node[sum,right=7mm of int] (s1) {$+$};
\node[block,right=7mm of s1] (Q) {$Q(\cdot)$};
\node[sum,right=8mm of Q] (mult) {$\times$};
\node[above=4mm of s1] (Sin) {$S(t)$};
\node[above=4mm of mult] (Min) {$M(t)$};
\node[block] (zoh) at (0,-14mm) {ZOH};
\node[block,align=center,anchor=west,font=\scriptsize] (etm) at ($(int.west)+(0,-14mm)$) {Event generator (\ref{eq:tk})\\ event iff $|e(t)|>\sigma|\hat{G}(t)|+\varepsilon$};
\draw[->] (K) -- node[above]{$u(t)$} (int);
\draw[->] (int) -- node[above]{$\hat{\theta}(t)$} (s1);
\draw[->] (s1) -- node[above]{$\theta(t)$} (Q);
\draw[->] (Q) -- node[above]{$y(t)$} (mult);
\draw[->] (Sin) -- (s1);
\draw[->] (Min) -- (mult);
\coordinate (out) at ($(mult.east)+(7mm,0)$);
\draw[->] (mult) -- node[above]{$\hat{G}(t)$} (out) |- ($(etm.south)+(0,-5mm)$) -| (zoh.south);
\draw[->] (out |- etm.east) -- (etm.east);
\draw[->,dashed] (etm.west) -- node[above,font=\scriptsize]{trigger} (zoh.east);
\draw[->] (zoh.north) -- node[right]{$\hat{G}(t_{k})$} (K.south);
\end{tikzpicture}
\caption{Robust event-triggered extremum seeking for a scalar map. The event generator monitors the measured gradient estimate $\hat{G}(t)$ and releases a new sample $\hat{G}(t_{k})$ to the zero-order hold (ZOH) only when the triggering condition is violated.}
\label{fig:BD}
\end{figure}
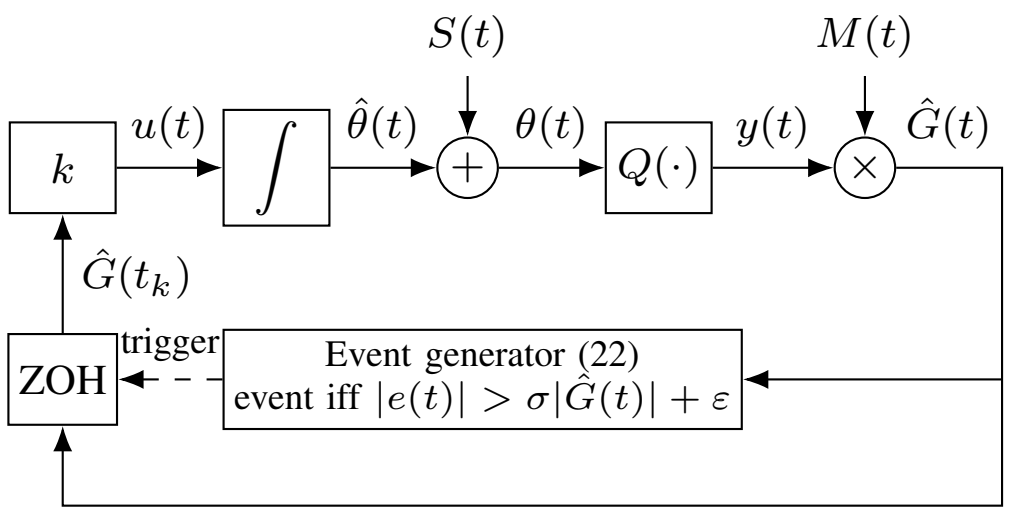

The scheme is depicted in Fig.~\ref{fig:BD}. A single sinusoid is used both to perturb the input and to demodulate the output \cite{K:2014}:
\begin{align}
S(t)=M(t)=a\sin(\omega t)\,, \quad a>0\,, \quad \omega>0\,, \label{eq:S}
\end{align}
with period
\begin{align}
T:=\frac{2\pi}{\omega}\,. \label{eq:T}
\end{align}
The dither is superimposed on the current estimate of the optimizer, $\theta(t)=\hat{\theta}(t)+S(t)$, where $\hat{\theta}(t)$ is the state of an integrator. With the estimation error
\begin{align}
\tilde{\theta}(t):=\hat{\theta}(t)-\theta^{\ast}\,, \label{eq:tildetheta}
\end{align}
one has $\theta(t)-\theta^{\ast}=\tilde{\theta}(t)+S(t)$ and, from (\ref{eq:map}),
\begin{align}
y(t)=Q^{\ast}+\frac{H^{\ast}}{2}\tilde{\theta}^{2}(t)+H^{\ast}S(t)\tilde{\theta}(t)+\frac{H^{\ast}}{2}S^{2}(t)\,. \label{eq:y_expanded}
\end{align}

\subsection{Structure of the Gradient Estimate}

The gradient estimate is obtained by demodulation,
\begin{align}
\hat{G}(t):=M(t)y(t)\,. \label{eq:Ghat_def}
\end{align}
Multiplying (\ref{eq:y_expanded}) by $M(t)=S(t)=a\sin(\omega t)$ and using the identities $2\sin^{2}(\omega t)=1-\cos(2\omega t)$ and $4\sin^{3}(\omega t)=3\sin(\omega t)-\sin(3\omega t)$, we arrive at the exact decomposition
\begin{align}
\hat{G}(t)=\frac{a^{2}}{2}H^{\ast}\beta(t)\,\tilde{\theta}(t)+d_{q}(t,\tilde{\theta}(t))+\Delta(t)\,, \label{eq:Ghat}
\end{align}
where
\begin{align}
\beta(t)&:=1-\cos(2\omega t)=2\sin^{2}(\omega t)\in[0,2]\,, \label{eq:beta}\\
d_{q}(t,\tilde{\theta})&:=\frac{a}{2}H^{\ast}\sin(\omega t)\,\tilde{\theta}^{2}\,, \label{eq:dq}\\
\Delta(t)&:=S(t)\left(Q^{\ast}+\frac{H^{\ast}}{2}S^{2}(t)\right) \nonumber\\
&\ =a\left(Q^{\ast}+\frac{3}{8}H^{\ast}a^{2}\right)\sin(\omega t)-\frac{H^{\ast}a^{3}}{8}\sin(3\omega t)\,. \label{eq:Delta}
\end{align}
The first term in (\ref{eq:Ghat}) is the useful part of the estimate: it is linear in $\tilde{\theta}$ and its coefficient has the sign of $H^{\ast}$ and average value $\frac{a^{2}}{2}H^{\ast}$ over one period. The term $d_{q}$ is quadratic in the error and vanishes at $\tilde{\theta}=0$. The term $\Delta(t)$ does not depend on $\tilde{\theta}$ and does not vanish: it is a persistent oscillation generated by the dither and by the offset $Q^{\ast}$ of the map. From (\ref{eq:Delta}),
\begin{align}
|\Delta(t)|\leq\bar{d}:=a\left(|Q^{\ast}|+\frac{a^{2}h}{2}\right)\,, \quad \forall t\geq 0\,. \label{eq:dbar}
\end{align}
Moreover, since $|\beta(t)|\leq 2$ and $|\sin(\omega t)|\leq 1$,
\begin{align}
|\hat{G}(t)|&\leq\ell\,|\tilde{\theta}(t)|+\frac{ah}{2}\tilde{\theta}^{2}(t)+\bar{d}\,, \quad \ell:=a^{2}h\,, \label{eq:Ghat_bound}\\
|d_{q}(t,\tilde{\theta})|&\leq\frac{ah}{2}\tilde{\theta}^{2}\,. \label{eq:dq_bound}
\end{align}

The next lemma gathers the three properties of the probing signal on which the whole analysis rests.

\begin{lemma}\label{lem:dither}
For all $t\geq 0$, the following statements hold.
\begin{enumerate}
\item[(i)] (Sign structure) $\beta(t)\geq 0$ and $\beta(t)\leq 2$.
\item[(ii)] (Persistence of excitation) $\displaystyle\int_{t}^{t+T}\beta(\tau)\,d\tau=T$.
\item[(iii)] (Bounded primitive) The $T$-periodic zero-mean primitive of $\Delta$, namely
\begin{align}
\Psi_{\Delta}(t):=&-\frac{a}{\omega}\left(Q^{\ast}+\frac{3}{8}H^{\ast}a^{2}\right)\cos(\omega t) \nonumber\\
&+\frac{H^{\ast}a^{3}}{24\,\omega}\cos(3\omega t)\,, \label{eq:Psi}
\end{align}
satisfies $\dot{\Psi}_{\Delta}(t)=\Delta(t)$ and
\begin{align}
|\Psi_{\Delta}(t)|\leq\frac{\kappa_{\Delta}}{\omega}\,, \quad \kappa_{\Delta}:=a|Q^{\ast}|+\frac{5}{12}ha^{3}\,. \label{eq:Psi_bound}
\end{align}
\end{enumerate}
\end{lemma}

\textit{Proof:} Items (i) and (ii) follow from $\beta(t)=2\sin^{2}(\omega t)$ and from $\int_{t}^{t+T}\cos(2\omega\tau)d\tau=0$, since $2\omega T=4\pi$. Item (iii) is obtained by integrating each sinusoid in (\ref{eq:Delta}) termwise and bounding $|Q^{\ast}+\frac{3}{8}H^{\ast}a^{2}|\leq|Q^{\ast}|+\frac{3}{8}ha^{2}$, so that $|\Psi_{\Delta}|\leq\frac{a}{\omega}(|Q^{\ast}|+\frac{3}{8}ha^{2}+\frac{1}{24}ha^{2})$. \hfill $\blacksquare$

Items (i) and (ii) state that the instantaneous adaptation gain never has the wrong sign and has a strictly positive mean; together they replace the averaging step of the classical analysis.

\subsection{Sample-and-Hold Emulation}

The adaptation gain is a scalar $k\in\R$ with
\begin{align}
\sgn(k)=-\sgn(H^{\ast})\,, \quad \bar{k}:=|k|\,, \label{eq:k}
\end{align}
so that $kH^{\ast}=-\bar{k}h$: a negative gain is used for minimization and a positive one for maximization. The estimator is the integrator $\dot{\hat{\theta}}(t)=u(t)$ and, since $\theta^{\ast}$ is constant, $\dot{\tilde{\theta}}(t)=u(t)$. The continuous-time law $u=k\hat{G}$ is emulated with a zero-order hold: the control value is recomputed only at the instants $(t_{k})_{k\in\N}$, $t_{0}=0$, delivered by the event generator,
\begin{align}
u(t)=k\,\hat{G}(t_{k})\,, \quad t\in[t_{k},t_{k+1})\,, \quad k\in\N\,. \label{eq:u}
\end{align}
The mismatch between the held and the current estimates is the sampling error
\begin{align}
e(t):=\hat{G}(t_{k})-\hat{G}(t)\,, \quad t\in[t_{k},t_{k+1})\,, \quad k\in\N\,. \label{eq:e}
\end{align}
Collecting (\ref{eq:Ghat}), (\ref{eq:u}) and (\ref{eq:e}), the estimation error obeys, on each interval $[t_{k},t_{k+1})$,
\begin{align}
\dot{\tilde{\theta}}(t)=-\frac{\bar{k}ha^{2}}{2}\beta(t)\,\tilde{\theta}(t)+k\,\Delta(t)+k\,\chi(t)\,, \label{eq:closed_loop}
\end{align}
where the lumped input
\begin{align}
\chi(t):=e(t)+d_{q}(t,\tilde{\theta}(t)) \label{eq:chi}
\end{align}
collects the sampling error, which is shaped by the event generator, and the quadratic term, which vanishes at the origin. Equation (\ref{eq:closed_loop}) is exact. Observe that $d_{q}$ appears twice: inside $\chi$, where it drives the error, and inside $\hat{G}$, which is the signal the event generator watches.

\subsection{Triggering Rule with Offset}

Because $\Delta(t)$ does not vanish at $\tilde{\theta}=0$, neither does $\hat{G}(t)$; in fact $\hat{G}(t)$ crosses zero at every half period. A purely relative rule of the form $|e|\leq\sigma|\hat{G}|$ would therefore force events to accumulate near those zero crossings. This is avoided by adding a positive offset, in the spirit of \cite[Sec.~3.4.1]{LZJ:2020}.

\begin{definition}[Static event generator with offset]\label{def:trigger}
Given $\sigma>0$ and $\varepsilon>0$, the sequence of event times $\{t_{k}\}_{k\in\N}$ is defined by $t_{0}=0$ and
\begin{align}
t_{k+1}:=\inf\left\{t>t_{k}:\ |e(t)|>\sigma|\hat{G}(t)|+\varepsilon\right\}\,, \label{eq:tk}
\end{align}
with the convention $t_{k+1}=+\infty$ when the set is empty. At each $t_{k}$ the control value (\ref{eq:u}) is refreshed.
\end{definition}

Implementing (\ref{eq:tk}) requires only $\hat{G}(t)$ and its last transmitted value $\hat{G}(t_{k})$; $\theta^{\ast}$, $Q^{\ast}$ and $H^{\ast}$ play no role in it. By construction, along the solutions,
\begin{align}
|e(t)|\leq\sigma|\hat{G}(t)|+\varepsilon\,, \quad \forall t\geq 0\,. \label{eq:e_bound}
\end{align}

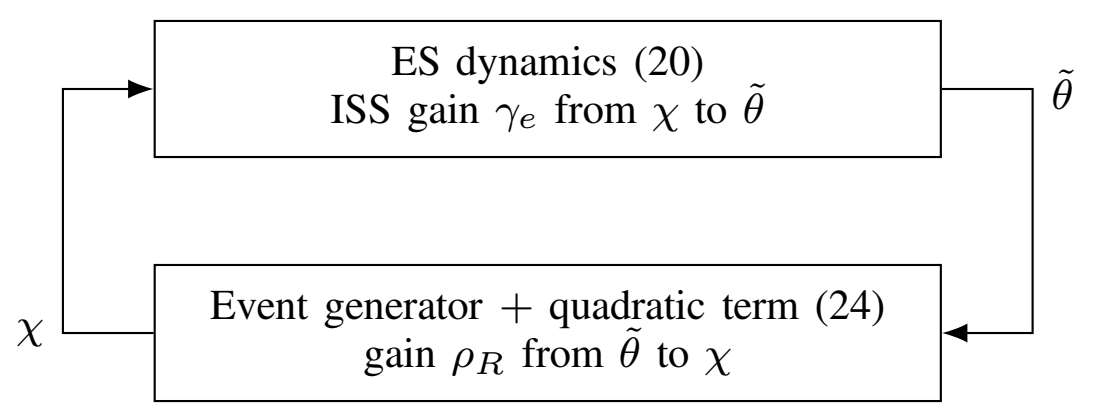
\begin{figure}[!t]
\centering
\begin{tikzpicture}[>=Latex,font=\footnotesize,block/.style={draw,align=center,minimum width=52mm,minimum height=9mm}]
\node[block] (A) {ES dynamics (\ref{eq:closed_loop})\\ ISS gain $\gamma_{e}$ from $\chi$ to $\tilde{\theta}$};
\node[block,below=7mm of A] (B) {Event generator $+$ quadratic term (\ref{eq:rho})\\ gain $\rho_{R}$ from $\tilde{\theta}$ to $\chi$};
\draw[->] (A.east) -- ++(6mm,0) node[right]{$\tilde{\theta}$} |- (B.east);
\draw[->] (B.west) -- ++(-6mm,0) node[left]{$\chi$} |- (A.west);
\end{tikzpicture}
\caption{Interconnection view of the closed loop. Stability holds under the small-gain condition $\gamma_{e}\rho_{R}<1$.}
\label{fig:SG}
\end{figure}

\textit{The event generator as a gain.} Combining (\ref{eq:e_bound}) with (\ref{eq:Ghat_bound}) and (\ref{eq:dq_bound}), the input (\ref{eq:chi}) satisfies
\begin{align}
|\chi(t)|\leq\sigma\ell\,|\tilde{\theta}(t)|+\frac{(1+\sigma)ah}{2}\tilde{\theta}^{2}(t)+\delta\,, \quad \delta:=\sigma\bar{d}+\varepsilon\,, \nonumber
\end{align}
the factor $1+\sigma$ being the price of the twofold appearance of $d_{q}$ noted above. Consequently, for any $R>0$, as long as $|\tilde{\theta}(t)|\leq R$,
\begin{align}
|\chi(t)|\leq\rho_{R}\,|\tilde{\theta}(t)|+\delta\,, \quad \rho_{R}:=\sigma a^{2}h+\frac{(1+\sigma)ah}{2}R\,. \label{eq:rho}
\end{align}
Following \cite{LJ:2015,LZJ:2020}, the event generator together with the quadratic term is viewed as an operator from $\tilde{\theta}$ to $\chi$ with linear gain $\rho_{R}$ and offset $\delta$ on the ball of radius $R$. Feeding $\chi$ into (\ref{eq:closed_loop}), which returns $\tilde{\theta}$, closes the loop of Fig.~\ref{fig:SG}. Section~IV quantifies the gain of each block and shows that, when their product is small enough, the interconnection is stable.

\section{MAIN RESULTS}

Throughout, we use the excitation level
\begin{align}
\lambda:=a^{2}h>0\,, \label{eq:lambda}
\end{align}
which is the mean value over one period of $a^{2}h\beta(t)$; equivalently, $\bar{k}\lambda/2$ is the average damping coefficient of (\ref{eq:closed_loop}).

\subsection{ISS of the Extremum Seeking Dynamics}

We first treat $\chi$ in (\ref{eq:closed_loop}) as an arbitrary exogenous signal, disregarding how it is produced. In the framework of \cite{LJ:2015,LZJ:2020}, the resulting estimate is the property that the plant is assumed to possess; here it is derived rather than assumed. Item (a) below is an ISS-Lyapunov inequality for the discrete-time system obtained by sampling the Lyapunov function once per dither period, in the sense of \cite{JW:2002}.

\begin{lemma}[ISS with respect to $\chi$]\label{lem:ISS}
Let Assumption~\ref{ass:map} hold and let the gain satisfy (\ref{eq:k}). Consider (\ref{eq:closed_loop}) driven by an arbitrary piecewise continuous input $\chi$, and define
\begin{align}
z(t)&:=\tilde{\theta}(t)-k\,\Psi_{\Delta}(t)\,, \quad W(t):=z^{2}(t)\,, \label{eq:z_W}\\
m_{0}&:=\frac{\bar{k}\lambda}{8}\,, \quad \gamma_{e}:=\frac{48}{\lambda}\,, \quad \delta_{\omega}:=\frac{\bar{k}\kappa_{\Delta}}{\omega}\,, \quad \bar{\nu}:=\bar{k}^{2}\lambda\frac{\kappa_{\Delta}}{\omega}\,. \label{eq:constants}
\end{align}
If
\begin{align}
\omega\geq\omega_{0}:=16\pi\,\bar{k}\lambda\,, \label{eq:omega0}
\end{align}
then, for all $t\geq 0$:
\begin{enumerate}
\item[(a)] (One-period ISS-Lyapunov inequality)
\begin{align}
W(t+T)\leq\left(1-m_{0}T\right)W(t)+T\,\frac{64}{\bar{k}\lambda}\,\eta_{t}^{2}\,, \label{eq:ISS_W}
\end{align}
with $\eta_{t}:=\bar{k}\sup_{[t,t+T]}|\chi|+\bar{\nu}$;
\item[(b)] (ISS estimate)
\begin{align}
|\tilde{\theta}(t)|\leq 4\exp\left(-\frac{m_{0}}{2}t\right)|\tilde{\theta}(0)|+\gamma_{e}\sup_{[0,t]}|\chi|+53\,\delta_{\omega}\,. \label{eq:ISS_gain}
\end{align}
Moreover, the right-hand side of (\ref{eq:ISS_gain}) involves $\chi$ only through its restriction to $[0,t]$.
\end{enumerate}
\end{lemma}

\textit{Proof:} \textit{Step 1 (Removing the oscillatory forcing).} Although $k\Delta(t)$ never dies out, it integrates to zero over each period. Since $\dot{\Psi}_{\Delta}=\Delta$, the variable $z$ in (\ref{eq:z_W}) satisfies
\begin{align}
\dot{z}(t)=-\frac{\bar{k}\lambda}{2}\beta(t)\,z(t)+k\chi(t)+\nu(t)\,, \label{eq:dotz}\\
\nu(t):=-\frac{k\bar{k}\lambda}{2}\beta(t)\Psi_{\Delta}(t)\,, \nonumber
\end{align}
where we used $-\frac{\bar{k}ha^{2}}{2}\beta\,\tilde{\theta}=-\frac{\bar{k}\lambda}{2}\beta(z+k\Psi_{\Delta})$. From Lemma~\ref{lem:dither}(i) and (\ref{eq:Psi_bound}),
\begin{align}
|\nu(t)|\leq\bar{\nu}\,, \quad |\tilde{\theta}(t)-z(t)|\leq\delta_{\omega}\,. \label{eq:nu_bound}
\end{align}

\textit{Step 2 (Derivative of $W$).} Along (\ref{eq:dotz}), for $\tau\in[t,t+T]$,
\begin{align}
\dot{W}(\tau)&=-\bar{k}\,q(\tau)+2z(\tau)\left(k\chi(\tau)+\nu(\tau)\right) \nonumber\\
&\leq-\bar{k}\,q(\tau)+2\eta_{t}|z(\tau)|\,, \label{eq:dotW}
\end{align}
with $q(\tau):=\lambda\beta(\tau)z^{2}(\tau)\geq 0$. Since $\beta$ touches zero twice per period, the first term cannot provide a pointwise decrease of $W$; the decrease emerges only after integration over a period, where Lemma~\ref{lem:dither}(ii) can be exploited.

\textit{Step 3 (Decrease over one period).} Let $L_{0}:=\bar{k}\lambda$. By (\ref{eq:T}) and (\ref{eq:omega0}), $L_{0}T=2\pi\bar{k}\lambda/\omega\leq 1/8$. From (\ref{eq:dotz}) and Lemma~\ref{lem:dither}(i), $|\dot{z}(\tau)|\leq L_{0}|z(\tau)|+\eta_{t}$ on $[t,t+T]$, so Gr\"onwall's inequality \cite{K:2002} and $e^{1/8}<2$ give
\begin{align}
|z(\tau)|&\leq 2\left(|z(t)|+\eta_{t}T\right)\,, \label{eq:z_window}\\
|z(\tau)-z(t)|&\leq 2L_{0}T|z(t)|+2\eta_{t}T\,. \label{eq:z_var}
\end{align}
For fixed $\tau$, $x\mapsto\sqrt{\lambda\beta(\tau)}\,|x|$ is a seminorm, hence $q(\tau)\geq\frac{1}{2}\lambda\beta(\tau)z^{2}(t)-\lambda\beta(\tau)\left(z(\tau)-z(t)\right)^{2}$. Integrating over $[t,t+T]$, using $\beta\leq 2$, Lemma~\ref{lem:dither}(ii), (\ref{eq:z_var}) and $L_{0}T\leq 1/8$,
\begin{align}
\int_{t}^{t+T}q(\tau)d\tau&\geq\frac{T\lambda}{2}z^{2}(t)-16\lambda T^{3}\left(L_{0}^{2}z^{2}(t)+\eta_{t}^{2}\right) \nonumber\\
&\geq\frac{T\lambda}{4}z^{2}(t)-16\lambda T^{3}\eta_{t}^{2}\,. \label{eq:q_int}
\end{align}
Also, from (\ref{eq:z_window}), $\int_{t}^{t+T}|z|d\tau\leq 2T(|z(t)|+\eta_{t}T)$. Integrating (\ref{eq:dotW}) and using $16\bar{k}\lambda T=16L_{0}T\leq 2$,
\begin{align}
W(t+T)-W(t)\leq-\frac{\bar{k}\lambda T}{4}z^{2}(t)+4\eta_{t}T|z(t)|+6\eta_{t}^{2}T^{2}\,. \label{eq:dW1}
\end{align}
By Young's inequality $xy\leq\frac{x^{2}}{2\varpi}+\frac{\varpi y^{2}}{2}$ with $x=4\eta_{t}$, $y=|z(t)|$ and $\varpi=\bar{k}\lambda/4$,
\begin{align}
W(t+T)\leq\left(1-\frac{\bar{k}\lambda T}{8}\right)W(t)+T\frac{32}{\bar{k}\lambda}\eta_{t}^{2}+6T^{2}\eta_{t}^{2}\,, \nonumber
\end{align}
and since $6T\leq 32/(\bar{k}\lambda)$ (because $L_{0}T\leq 1/8<16/3$), inequality (\ref{eq:ISS_W}) follows.

\textit{Step 4 (ISS estimate).} Fix $t>0$ and set $\bar{\chi}:=\sup_{[0,t]}|\chi|$, $\bar{\eta}:=\bar{k}\bar{\chi}+\bar{\nu}$ and $\mu_{0}:=64\bar{\eta}^{2}/(\bar{k}\lambda)$. Whenever $[iT,(i+1)T]\subset[0,t]$, the window constant obeys $\eta_{iT}\leq\bar{\eta}$, and (\ref{eq:ISS_W}) becomes the linear recursion $W((i+1)T)\leq cW(iT)+T\mu_{0}$, whose coefficient $c:=1-m_{0}T$ lies in $(0,1)$ because $m_{0}T=L_{0}T/8\leq 1/64$. Iterating, for all $j\in\N$ with $jT\leq t$,
\begin{align}
W(jT)\leq c^{j}W(0)+\frac{\mu_{0}}{m_{0}}\left(1-c^{j}\right)\leq e^{-m_{0}jT}W(0)+\frac{\mu_{0}}{m_{0}}\,, \label{eq:WjT}
\end{align}
where $1-x\leq e^{-x}$ was used. Now take $\tau\in[jT,\min\{(j+1)T,t\}]$. By (\ref{eq:z_window}) on $[jT,\tau]$, (\ref{eq:WjT}), $\sqrt{x+y}\leq\sqrt{x}+\sqrt{y}$ and $e^{m_{0}T/2}\leq 2$,
\begin{align}
|z(\tau)|\leq 4e^{-\frac{m_{0}}{2}\tau}|z(0)|+2\sqrt{\frac{\mu_{0}}{m_{0}}}+2\bar{\eta}T\,. \label{eq:z_final}
\end{align}
Here $2\sqrt{\mu_{0}/m_{0}}=2\sqrt{512}\,\bar{\eta}/(\bar{k}\lambda)\leq 46\,\bar{\eta}/(\bar{k}\lambda)$ and, from $L_{0}T\leq 1$, $2\bar{k}T\leq 2/\lambda$. Substituting $\bar{\eta}=\bar{k}\bar{\chi}+\bar{\nu}$, the coefficient of $\bar{\chi}$ in (\ref{eq:z_final}) is at most $48/\lambda=\gamma_{e}$, while the terms in $\bar{\nu}$ add up to at most $46\bar{\nu}/(\bar{k}\lambda)+2\bar{\nu}T\leq 48\,\delta_{\omega}$, by the definition of $\bar{\nu}$ in (\ref{eq:constants}) and $\bar{k}\lambda T\leq 1$. Finally, (\ref{eq:nu_bound}) gives $|\tilde{\theta}(\tau)|\leq|z(\tau)|+\delta_{\omega}$ and $|z(0)|\leq|\tilde{\theta}(0)|+\delta_{\omega}$, which contributes $5\delta_{\omega}$ and yields (\ref{eq:ISS_gain}). Only the values of $\chi$ on $[0,t]$ were used. \hfill $\blacksquare$

\subsection{Small-Gain Theorem}

\begin{theorem}\label{thm:main}
Let Assumption~\ref{ass:map} hold and consider the map (\ref{eq:map}) in feedback with the held gradient law (\ref{eq:u}), (\ref{eq:k}) and the event generator of Definition~\ref{def:trigger}. Suppose that the parameters $R>0$, $\sigma>0$ and $\varepsilon>0$ fulfill the small-gain condition
\begin{align}
\kappa:=\gamma_{e}\,\rho_{R}<1\,, \quad\text{i.e.,}\quad \sigma a^{2}h+\frac{(1+\sigma)ah}{2}R<\frac{a^{2}h}{48}\,, \label{eq:SG}
\end{align}
and the offset condition
\begin{align}
\gamma_{e}\,\delta\leq\frac{(1-\kappa)R}{4}\,, \quad\text{i.e.,}\quad \sigma\bar{d}+\varepsilon\leq\frac{(1-\kappa)a^{2}hR}{192}\,, \label{eq:offset}
\end{align}
with $\gamma_{e}$, $\rho_{R}$, $\delta$ and $\bar{d}$ in (\ref{eq:constants}), (\ref{eq:rho}) and (\ref{eq:dbar}). Let the dither frequency satisfy
\begin{align}
\omega>\omega^{\ast}:=\max\left\{16\pi\bar{k}\left(\lambda+\rho_{R}\right),\ \frac{212\,\bar{k}\kappa_{\Delta}}{(1-\kappa)R}\right\}\,. \label{eq:omega_star}
\end{align}
Then, for every initial condition with $|\theta(0)-\theta^{\ast}|\leq(1-\kappa)R/8$, the solution exists for all $t\geq 0$, satisfies $|\tilde{\theta}(t)|<R$, and
\begin{align}
|\theta(t)-\theta^{\ast}|&\leq M_{\theta}\,e^{-\frac{m}{2}t}+a+\bar{r}\,, \label{eq:theta_bound}\\
|y(t)-Q^{\ast}|&\leq M_{y}\,e^{-mt}+2h\left(a^{2}+\bar{r}^{2}\right)\,, \label{eq:y_bound}
\end{align}
where $m:=(1-\kappa^{2})\bar{k}\lambda/8$, $M_{\theta}:=4|\theta(0)-\theta^{\ast}|$, $M_{y}:=hM_{\theta}^{2}$ and
\begin{align}
\bar{r}\leq\frac{3\gamma_{e}\delta}{\sqrt{1-\kappa^{2}}}+\mathcal{O}\left(\frac{1}{\omega}\right)\,. \label{eq:rbar}
\end{align}
In addition, the inter-event times satisfy $t_{k+1}-t_{k}\geq\tau^{\ast}>0$ for all $k\in\N$, with $\tau^{\ast}$ given in (\ref{eq:tau_star}) below, so that Zeno behavior does not occur.
\end{theorem}

\textit{Proof:} Let $t_{\infty}:=\lim_{k\to\infty}t_{k}\leq+\infty$. Between events the control is constant, so $\tilde{\theta}$ is continuous and piecewise affine on $[0,t_{\infty})$ and cannot escape in finite time. Part~A shows that $|\tilde{\theta}(t)|<R$ on $[0,t_{\infty})$ and establishes the convergence estimates; Part~B uses only the bound $|\tilde{\theta}|<R$ to prove that $t_{\infty}=+\infty$.

\textit{A. Stability analysis.}

\textit{Step 1 (The ball $|\tilde{\theta}|<R$ is invariant).} Since $S(0)=0$, $\tilde{\theta}(0)=\theta(0)-\theta^{\ast}$ and thus $4|\tilde{\theta}(0)|\leq(1-\kappa)R/2$. By (\ref{eq:omega_star}), $\omega\geq\omega_{0}$ and $53\delta_{\omega}=53\bar{k}\kappa_{\Delta}/\omega<(1-\kappa)R/4$. If the claim were false, continuity of $\tilde{\theta}$ would provide a first time\linebreak $t_{e}\in(0,t_{\infty})$ at which $|\tilde{\theta}(t_{e})|=R$. Up to that time the bound (\ref{eq:rho}) is in force, giving $|\chi|\leq\rho_{R}R+\delta$ on $[0,t_{e}]$. Because (\ref{eq:ISS_gain}) only sees $\chi$ on $[0,t_{e}]$ and $\chi$ is piecewise continuous, we may evaluate Lemma~\ref{lem:ISS}(b) at $t=t_{e}$ (continuing $\chi$ by zero beyond $t_{e}$, for definiteness). With (\ref{eq:SG}) and (\ref{eq:offset}),
\begin{align}
|\tilde{\theta}(t_{e})|&\leq 4|\tilde{\theta}(0)|+\kappa R+\gamma_{e}\delta+53\delta_{\omega} \nonumber\\
&<\frac{(1-\kappa)R}{2}+\kappa R+\frac{(1-\kappa)R}{4}+\frac{(1-\kappa)R}{4}=R\,, \nonumber
\end{align}
a contradiction. Thus (\ref{eq:rho}) can be used freely on $[0,t_{\infty})$.

\textit{Step 2 (Closing the loop).} From (\ref{eq:rho}) and (\ref{eq:nu_bound}), $|\chi(\tau)|\leq\rho_{R}|z(\tau)|+\rho_{R}\delta_{\omega}+\delta$. Then (\ref{eq:dotz}) gives $|\dot{z}(\tau)|\leq L|z(\tau)|+\hat{\eta}$ with
\begin{align}
L:=\bar{k}(\lambda+\rho_{R})\,, \quad \hat{\eta}:=\bar{k}\left(\delta+\rho_{R}\delta_{\omega}\right)+\bar{\nu}\,. \nonumber
\end{align}
By (\ref{eq:omega_star}), $LT\leq 1/8$, so Gr\"onwall's inequality yields $|z(\tau)|\leq 2(|z(t)|+\hat{\eta}T)$ on $[t,t+T]$ and the input of (\ref{eq:ISS_W}) obeys
\begin{align}
\eta_{t}\leq 2\bar{k}\rho_{R}|z(t)|+\eta'\,, \quad \eta':=\bar{k}\delta+\bar{k}\rho_{R}\left(2\hat{\eta}T+\delta_{\omega}\right)+\bar{\nu}\,. \label{eq:eta_prime}
\end{align}

\textit{Step 3 (Small gain).} Inserting (\ref{eq:eta_prime}) into (\ref{eq:ISS_W}) and using the weighted Young inequality $(x+y)^{2}\leq\frac{9}{8}x^{2}+9y^{2}$,
\begin{align}
W(t+T)\leq\left(1-m_{0}T+T\frac{288\,\bar{k}\rho_{R}^{2}}{\lambda}\right)W(t)+T\frac{576\,\eta'^{2}}{\bar{k}\lambda}\,. \label{eq:W_SG}
\end{align}
The quantity $288\bar{k}\rho_{R}^{2}/\lambda$ is the loop gain of the interconnection in Fig.~\ref{fig:SG} as seen by the Lyapunov function. Since $\rho_{R}=\kappa\lambda/48$ by (\ref{eq:SG}), this loop gain equals $\kappa^{2}\bar{k}\lambda/8=\kappa^{2}m_{0}$, and $\kappa<1$ is exactly what is needed to retain a positive decrease rate. Therefore, with
\begin{align}
m:=(1-\kappa^{2})\,m_{0}=\frac{(1-\kappa^{2})\bar{k}\lambda}{8}\,, \quad \mu:=\frac{576\,\eta'^{2}}{\bar{k}\lambda}\,, \label{eq:m_mu}
\end{align}
we obtain $W(t+T)\leq(1-mT)W(t)+T\mu$, with $mT\leq m_{0}T\leq 1/64$.

\textit{Step 4 (Convergence).} The recursion $W(t+T)\leq(1-mT)W(t)+T\mu$ has the same form as in Step~4 of Lemma~\ref{lem:ISS}; running that argument again with $(m,\mu,\hat{\eta})$ instead of $(m_{0},\mu_{0},\bar{\eta})$ gives
\begin{align}
|\tilde{\theta}(t)|&\leq 4e^{-\frac{m}{2}t}|\tilde{\theta}(0)|+\bar{r}\,, \label{eq:tildetheta_final}\\
\bar{r}&:=\frac{136\,\eta'}{\bar{k}\lambda\sqrt{1-\kappa^{2}}}+2\hat{\eta}T+5\delta_{\omega}\,, \nonumber
\end{align}
where $2\sqrt{\mu/m}\!=\!2\sqrt{4608}\,\eta'/(\bar{k}\lambda\sqrt{1-\kappa^{2}})\!\leq\! 136\,\eta'/(\bar{k}\lambda\sqrt{1-\kappa^{2}})$ was used. All three ingredients of $\bar{r}$ are explicit: $\eta'=\bar{k}\delta+\mathcal{O}(1/\omega)$, $\hat{\eta}T=\mathcal{O}(1/\omega)$ and $\delta_{\omega}=\mathcal{O}(1/\omega)$; since $136/\lambda\leq 3\gamma_{e}$, this is (\ref{eq:rbar}). 

\textit{Step 5 (Bounds on $\theta$ and $y$).} From $\theta-\theta^{\ast}=\tilde{\theta}+S$, $|S|\leq a$ and (\ref{eq:tildetheta_final}), estimate (\ref{eq:theta_bound}) holds with $M_{\theta}=4|\theta(0)-\theta^{\ast}|$. Moreover, $|y-Q^{\ast}|=\frac{h}{2}(\theta-\theta^{\ast})^{2}\leq\frac{h}{2}\left(M_{\theta}e^{-\frac{m}{2}t}+\bar{r}+a\right)^{2}$, and two applications of $(x+y)^{2}\leq 2x^{2}+2y^{2}$ give (\ref{eq:y_bound}) with $M_{y}=hM_{\theta}^{2}$.\\

\textit{B. Exclusion of Zeno behavior.}

We now use only the invariance established in Step~1. With $\phi(t):=\tilde{\theta}(t)+S(t)=\theta(t)-\theta^{\ast}$ one has $|\phi|\leq R+a$ and $y=Q^{\ast}+\frac{H^{\ast}}{2}\phi^{2}$, whence\\
\begin{align}
|y(t)|\leq\bar{y}:=|Q^{\ast}|+\frac{h}{2}(R+a)^{2}\,, \quad |\hat{G}(t)|\leq a\bar{y}\,. \nonumber
\end{align}
Writing $S_{k}:=S(t_{k})$, $y_{k}:=y(t_{k})$, $\phi_{k}:=\phi(t_{k})$, the sampling error (\ref{eq:e}) splits as
\begin{align}
e(t)=\left(S_{k}-S(t)\right)y_{k}+S(t)\left(y_{k}-y(t)\right)\,, \quad t\in[t_{k},t_{k+1})\,. \nonumber
\end{align}
Since $|S_{k}-S(t)|\leq a\omega(t-t_{k})$, $|y_{k}-y(t)|=\frac{h}{2}|\phi_{k}-\phi(t)||\phi_{k}+\phi(t)|\leq h(R+a)|\phi_{k}-\phi(t)|$ and, because $\dot{\tilde{\theta}}=k\hat{G}(t_{k})$ is constant between events, $|\phi_{k}-\phi(t)|\leq(\bar{k}a\bar{y}+a\omega)(t-t_{k})$, we conclude that
\begin{align}
|e(t)|&\leq L_{e}\,(t-t_{k})\,, \quad t\in[t_{k},t_{k+1})\,, \label{eq:e_linear}\\
L_{e}&:=a\omega\bar{y}+ah(R+a)\left(\bar{k}a\bar{y}+a\omega\right)\,. \nonumber
\end{align}
At an event, by (\ref{eq:tk}) and continuity of $e$ on $[t_{k},t_{k+1})$, $|e(t_{k+1}^{-})|\geq\sigma|\hat{G}(t_{k+1})|+\varepsilon\geq\varepsilon$, while $e(t_{k})=0$. Hence (\ref{eq:e_linear}) gives
\begin{align}
t_{k+1}-t_{k}\geq\tau^{\ast}:=\frac{\varepsilon}{L_{e}}>0\,, \quad \forall k\in\N\,. \label{eq:tau_star}
\end{align}
As $\tau^{\ast}$ does not depend on $k$, either finitely many events occur or $t_{k}\to\infty$; in both cases $t_{\infty}=+\infty$, the solution is defined for all $t\geq 0$ and Steps~2--5 of Part~A hold globally in time. \hfill $\blacksquare$

\subsection{Discussion}

\begin{remark}[Small-gain reading]
In (\ref{eq:SG}) two gains are multiplied: $\gamma_{e}$, which measures how strongly the ES dynamics amplify $\chi$ into $\tilde{\theta}$ (Lemma~\ref{lem:ISS}), and $\rho_{R}$, which measures how strongly the event generator and the quadratic term send $\tilde{\theta}$ back into $\chi$ (see (\ref{eq:rho})). This is the one-dimensional, linear-gain version of the composition condition $\underline{\alpha}^{-1}\circ\gamma\circ\rho<\mathrm{Id}$ of \cite{LJ:2015,LZJ:2020}. Two parameters feed $\rho_{R}$: the relative tolerance $\sigma$, and the size $R$ of the region where the result is claimed, which enters because the quadratic part of the map is treated as an input. It is worth observing that $\bar{k}$ drops out of $\gamma_{e}=48/(a^{2}h)$, whereas a larger excitation level $\lambda$ lowers it. Condition (\ref{eq:SG}) is in the standard form $\kappa<1$, and the strict inequality is used only to keep the decrease rate $m=(1-\kappa^{2})m_{0}$ positive; accordingly, the decay rate, the admissible offset, the admissible initial error and the residual radius in (\ref{eq:rbar}) all degrade as $\kappa\to 1$, the factor $(1-\kappa^{2})^{-1/2}$ in (\ref{eq:rbar}) being the counterpart of the familiar $(1-\kappa)^{-1}$ of the classical small-gain estimate $\gamma_{e}\delta/(1-\kappa)$. Had the cross term in Step~3 of the proof been split with the symmetric inequality $(x+y)^{2}\leq 2x^{2}+2y^{2}$, the loop gain would have read $512\bar{k}\rho_{R}^{2}/\lambda$ and only $\kappa<3/4$ could have been certified; the weighted splitting removes this artificial restriction at the cost of a larger constant multiplying $\eta'^{2}$.
\end{remark}

\begin{remark}[Design procedure]\label{rmk:design}
Dividing (\ref{eq:SG}) by $ah$, the small-gain condition becomes $\sigma a+\frac{1+\sigma}{2}R<\frac{a}{48}$, which suggests the following order of choices:
\begin{enumerate}
\item the gain $k$ and the amplitude $a$, which fix $\lambda$, $\ell$, $\bar{d}$ and $\gamma_{e}$;
\item a provisional threshold $\sigma_{0}<1/48$;
\item a radius
\begin{align}
0<R<\frac{2a\left(1-48\sigma_{0}\right)}{48\left(1+\sigma_{0}\right)}\,, \label{eq:R_choice}
\end{align}
which reserves an explicit share of the small-gain margin for the quadratic term, and the corresponding $\kappa_{0}:=\gamma_{e}\rho_{R}|_{\sigma=\sigma_{0}}<1$;
\item the threshold $\sigma\leq\min\{\sigma_{0},\,(1-\kappa_{0})a^{2}hR/(384\,\bar{d})\}$ and the offset $\varepsilon\leq(1-\kappa_{0})a^{2}hR/384$, so that (\ref{eq:offset}) holds;
\item the frequency $\omega>\omega^{\ast}$ in (\ref{eq:omega_star}).
\end{enumerate}
Since $\rho_{R}$, and hence $\kappa$, grows with $\sigma$, lowering $\sigma$ in step 4) preserves (\ref{eq:SG}) and gives $\kappa\leq\kappa_{0}$. Inequality (\ref{eq:R_choice}) forces $R<a/24$, i.e., the region in which convergence is certified is comparable to the probing amplitude. This is inherent to a statement proved on the exact, non-averaged loop. Being sufficient conditions obtained through worst-case bounds, (\ref{eq:SG})--(\ref{eq:omega_star}) are conservative, and Section~\ref{sec:sim} shows satisfactory behavior for parameters well beyond them.
\end{remark}

\begin{remark}[Trade-offs]\label{rmk:tradeoffs}
Three mechanisms determine the size of the final neighborhood in (\ref{eq:theta_bound}): the probing amplitude $a$, which is unavoidable in ES; the oscillation $\Delta(t)$, whose effect scales as $1/\omega$; and the event generator, which contributes a term proportional to $\gamma_{e}\delta$. Without a strictly positive $\varepsilon$ no uniform bound of the type (\ref{eq:tau_star}) could hold, because $\hat{G}$ keeps oscillating through zero even when $\tilde{\theta}=0$; what is obtained is practical stability in the presence of a non-vanishing matched disturbance, a situation also addressed in \cite{KK:2024}. Since $L_{e}$ in (\ref{eq:e_linear}) grows linearly with $\omega$, one has $\tau^{\ast}=\mathcal{O}(\varepsilon/\omega)$: a faster dither shrinks the residual set but demands faster updates, in agreement with the intuition that the held control must track the probing oscillation.
\end{remark}

\section{NUMERICAL EXAMPLE}\label{sec:sim}

The scalar map (\ref{eq:map}) is simulated with
\begin{align}
H^{\ast}=4\,, \quad Q^{\ast}=1\,, \quad \theta^{\ast}=1\,, \nonumber
\end{align}
so that $Q$ has a minimum and, by (\ref{eq:k}), a negative gain must be used. The probing signal (\ref{eq:S}) has amplitude $a=0.2$ and frequency $\omega=3$~rad/s, the adaptation gain is $k=-1$, and the integrator starts at $\hat{\theta}(0)=2$, i.e., one unit away from the optimizer. The event generator of Definition~\ref{def:trigger} is run with $\sigma=0.1$ and $\varepsilon=0.01$, and the triggering condition (\ref{eq:tk}) is evaluated at every step of a fixed-step integration with step size $10^{-3}$~s over a horizon of $200$~s. For these values the excitation level is $\lambda=a^{2}h=0.16$, so the mean time constant of (\ref{eq:closed_loop}) is $2/(\bar{k}\lambda)=12.5$~s; moreover, $\bar{d}=0.216$ and $\delta=\sigma\bar{d}+\varepsilon\approx 0.032$. The relative threshold $\sigma=0.1$ lies well above the value $1/48$ allowed by Remark~\ref{rmk:design}, so the example also serves to show how conservative the sufficient conditions are.

\begin{figure}[!t]
\centering
\includegraphics[width=0.82\columnwidth]{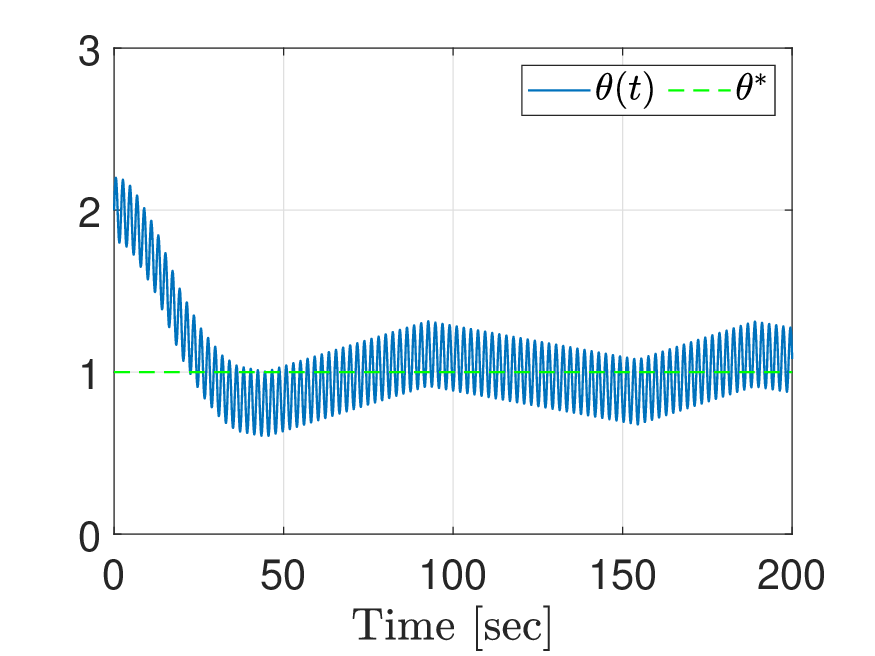}
\caption{Input $\theta(t)$ of the map and optimizer $\theta^{\ast}$.}
\label{fig:sim_theta}
\end{figure}

\begin{figure}[!t]
\centering
\includegraphics[width=0.82\columnwidth]{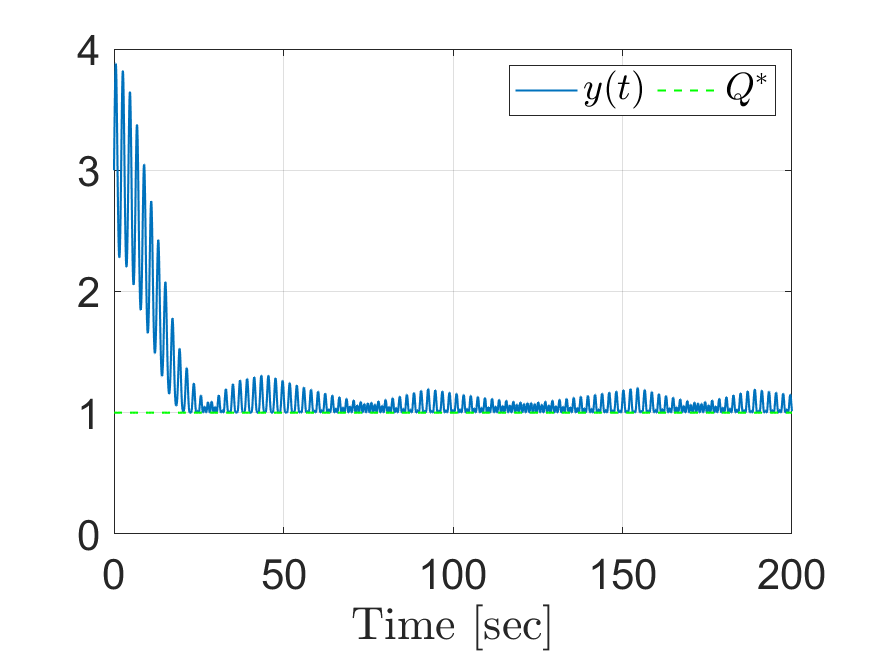}
\caption{Output $y(t)$ of the map and minimum value $Q^{\ast}$.}
\label{fig:sim_y}
\end{figure}

Figs.~\ref{fig:sim_theta} and \ref{fig:sim_y} show the input and the output of the map. Starting from $y(0)=3$, the input $\theta(t)$ descends toward $\theta^{\ast}$ within roughly $25$~s, a duration compatible with the time constant dictated by $\lambda$, and the output is brought from its initial excursions, which reach almost $4$ because of the dither, to the vicinity of $Q^{\ast}=1$. After the transient, $\theta(t)$ carries the fast oscillation of amplitude $a=0.2$ imposed by the probing signal, which is the irreducible term $a$ in (\ref{eq:theta_bound}), superimposed on a slow, low-amplitude excursion about $\theta^{\ast}$: the mean value of $\theta$ undershoots to about $0.75$ near $t=45$~s, then drifts to about $1.2$ near $t=90$~s and oscillates with decreasing amplitude afterwards. This slow component is the signature of the sample-and-hold implementation: between two events the integrator is driven by a constant $u=k\hat{G}(t_{k})$ which, owing to the persistent term $\Delta(t_{k})$ in (\ref{eq:Ghat}), is not exactly zero even when $\tilde{\theta}(t_{k})$ is small, so that $\hat{\theta}$ moves linearly until the triggering rule detects the discrepancy and refreshes the control. The excursion stays within approximately $0.25$ of $\theta^{\ast}$, and the corresponding output remains below $1.3$ for all $t>25$~s, in agreement with the practical-stability nature of Theorem~\ref{thm:main}, whose residual set is shaped not only by $a$ but also by the event-generator term in (\ref{eq:rbar}).

\begin{figure}[!t]
\centering
\includegraphics[width=0.82\columnwidth]{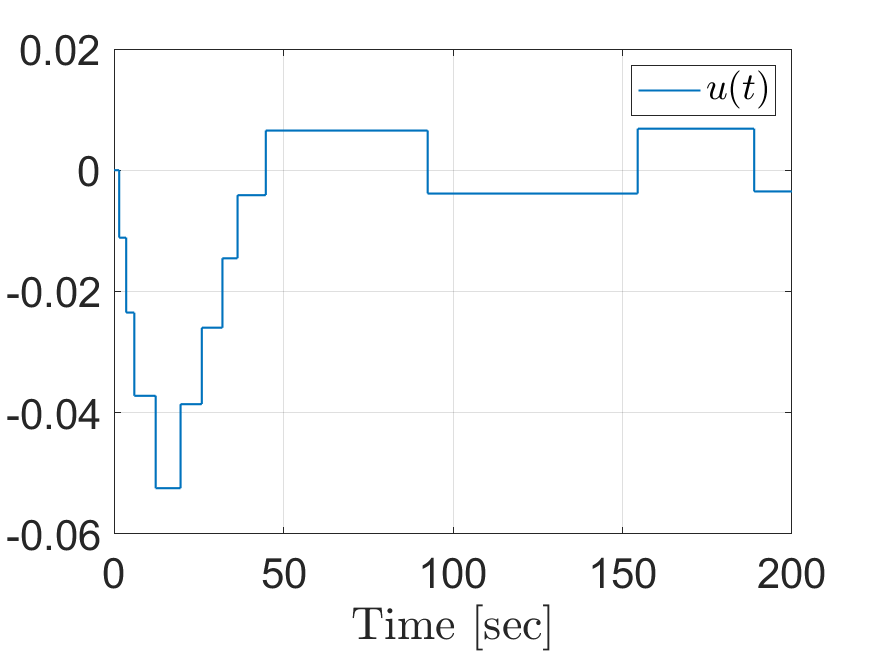}
\caption{Control signal $u(t)=k\hat{G}(t_{k})$ delivered by the zero-order hold.}
\label{fig:sim_u}
\end{figure}

\begin{figure}[!t]
\centering
\includegraphics[width=0.82\columnwidth]{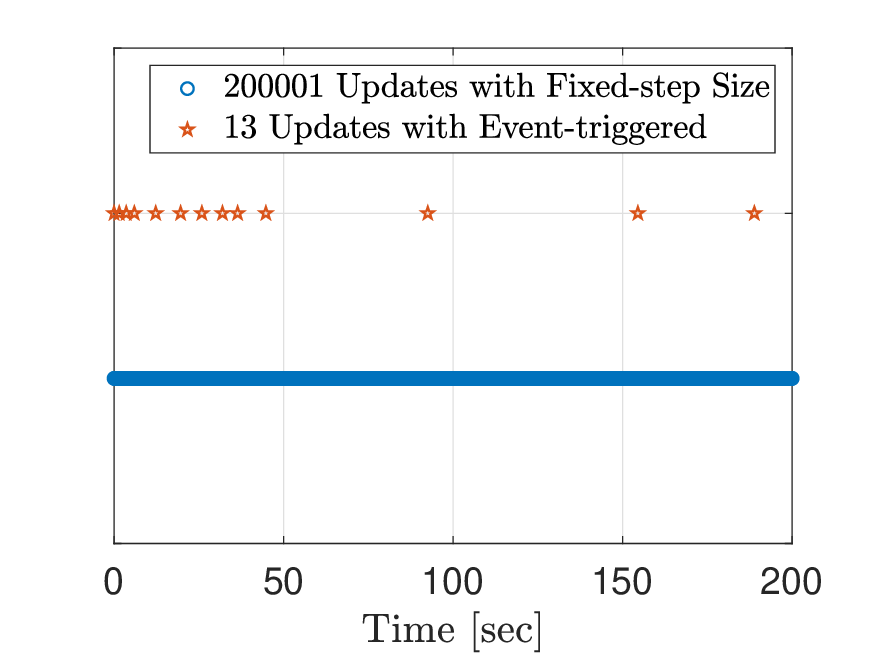}
\caption{Instants of control update: fixed-step implementation (circles, $200\,001$ updates) and event-triggered implementation (stars, $13$ updates).}
\label{fig:sim_updates}
\end{figure}

The control signal $u(t)=k\hat{G}(t_{k})$ is displayed in Fig.~\ref{fig:sim_u}. As imposed by (\ref{eq:u}), it is piecewise constant, with only $13$ distinct levels over the whole run. Its magnitude peaks at about $0.053$ around $t=15$~s, when the estimation error and hence the useful component $\frac{a^{2}}{2}H^{\ast}\beta(t)\tilde{\theta}$ of (\ref{eq:Ghat}) are largest, and then decays: after $t\approx 45$~s the held values do not exceed $0.008$ in absolute value. The distribution of the events, shown in Fig.~\ref{fig:sim_updates}, follows the same pattern. Ten of the thirteen updates occur during the transient, within the first $45$~s, while $|\hat{G}|$ is still appreciable and the relative term $\sigma|\hat{G}|$ in (\ref{eq:tk}) is easily exceeded; the remaining three occur near $t=92$, $155$ and $190$~s, separated by intervals of about one minute, during which the offset $\varepsilon$ dominates the threshold and keeps the hold active even though $\hat{G}(t)$ keeps oscillating through zero. No accumulation of events is observed anywhere, and the shortest interval between consecutive events is of the order of one second, orders of magnitude larger than the guaranteed bound (\ref{eq:tau_star}), which for $R=1.2$ (a ball that contains the whole trajectory, since $|\tilde{\theta}(0)|=1$) gives $\bar{y}=4.92$, $L_{e}\approx 4.7$ and $\tau^{\ast}=\varepsilon/L_{e}\approx 2$~ms.

Against the $200\,001$ updates of the time-triggered implementation, the event-triggered law converges to a neighborhood of the optimizer with $13$ updates, a reduction above $99.99\%$ in control transmissions. The price is the slow excursion of $\theta$ discussed above, which decreasing $\varepsilon$ and $\sigma$, or increasing $\omega$, would reduce at the expense of more updates (Remark~\ref{rmk:tradeoffs}).

\section{CONCLUSIONS}

We have analyzed an event-triggered implementation of gradient extremum seeking for scalar quadratic maps in which the triggering rule acts on the measured gradient estimate and includes a positive offset. The scalar structure of the demodulated signal, namely a nonnegative and persistently exciting adaptation coefficient plus an explicit zero-mean sinusoidal perturbation, allowed us to bypass averaging theory altogether: the ISS gain of the extremum seeking dynamics with respect to the sampling error was computed on the exact error equation, and closed-loop stability was obtained from a small-gain condition with the gain of the event generator, in which the quadratic term of the map is handled as part of the input. All design constants, the minimum dither frequency and the lower bound on the inter-event times are explicit. Dynamic triggers and Newton-based schemes are natural next steps.

\bibliographystyle{ieeetr}
\bibliography{refs}

@book{AK:2003,
  author    = {Ariyur, K. B. and Krsti{\'c}, M.},
  title     = {Real-Time Optimization by Extremum-Seeking Control},
  publisher = {Wiley},
  address   = {NJ, USA},
  year      = {2003}
}

@article{GKN:2012,
  author  = {Ghaffari, A. and Krsti{\'c}, M. and Ne{\v s}i{\'c}, D.},
  title   = {Multivariable {N}ewton-based extremum seeking},
  journal = {Automatica},
  volume  = {48},
  number  = {8},
  pages   = {1759--1767},
  year    = {2012}
}

@inproceedings{HJT:2012,
  author    = {Heemels, W. P. M. H. and Johansson, K. H. and Tabuada, P.},
  title     = {An introduction to event-triggered and self-triggered control},
  booktitle = {51st IEEE Conference on Decision and Control (CDC)},
  pages     = {3270--3285},
  year      = {2012}
}

@book{K:2002,
  author    = {Khalil, H. K.},
  title     = {Nonlinear Systems},
  edition   = {3rd},
  publisher = {Prentice Hall},
  address   = {Upper Saddle River, New Jersey},
  year      = {2002}
}

@incollection{K:2014,
  author    = {Krsti{\'c}, M.},
  title     = {Extremum seeking control},
  booktitle = {Encyclopedia of Systems and Control},
  editor    = {Baillieul, J. and Samad, T.},
  pages     = {413--416},
  publisher = {Springer},
  address   = {London},
  year      = {2014}
}

@article{KW:2000,
  author  = {Krsti{\'c}, M. and Wang, H.-H.},
  title   = {Stability of extremum seeking feedback for general nonlinear dynamic systems},
  journal = {Automatica},
  volume  = {36},
  pages   = {595--601},
  year    = {2000}
}

@article{LJ:2015,
  author  = {Liu, T. and Jiang, Z.-P.},
  title   = {A small-gain approach to robust event-triggered control of nonlinear systems},
  journal = {IEEE Transactions on Automatic Control},
  volume  = {60},
  number  = {8},
  pages   = {2072--2085},
  year    = {2015}
}

@book{LZJ:2020,
  author    = {Liu, T. and Zhang, P. and Jiang, Z.-P.},
  title     = {Robust Event-Triggered Control of Nonlinear Systems},
  publisher = {Springer},
  address   = {Singapore},
  year      = {2020}
}

@book{TRO_book2022,
  author    = {Oliveira, T. R. and Krsti{\'c}, M.},
  title     = {Extremum Seeking through Delays and {PDEs}},
  publisher = {Society for Industrial and Applied Mathematics},
  year      = {2022}
}

@article{OKT:2017,
  author  = {Oliveira, T. R. and Krsti{\'c}, M. and Tsubakino, D.},
  title   = {Extremum seeking for static maps with delays},
  journal = {IEEE Transactions on Automatic Control},
  volume  = {62},
  number  = {4},
  pages   = {1911--1926},
  year    = {2017}
}

@article{ROHDK:2025,
  author  = {Rodrigues, V. H. P. and Oliveira, T. R. and Hsu, L. and Diagne, M. and Krsti{\'c}, M.},
  title   = {Event-triggered and periodic event-triggered extremum seeking control},
  journal = {Automatica},
  volume  = {174},
  pages   = {112161},
  year    = {2025}
}

@article{T:2007,
  author  = {Tabuada, P.},
  title   = {Event-triggered real-time scheduling of stabilizing control tasks},
  journal = {IEEE Transactions on Automatic Control},
  volume  = {52},
  number  = {9},
  pages   = {1680--1685},
  year    = {2007}
}

@article{ZLJ:2017,
  author  = {Zhang, P. and Liu, T. and Jiang, Z.-P.},
  title   = {Input-to-state stabilization of nonlinear discrete-time systems with event-triggered controllers},
  journal = {Systems \& Control Letters},
  volume  = {103},
  pages   = {16--22},
  year    = {2017}
}

@article{ZAGSK:2007,
  author  = {Zhang, C. and Arnold, D. and Ghods, N. and Siranosian, A. and Krsti{\'c}, M.},
  title   = {Source seeking with non-holonomic unicycle without position measurement and with tuning of forward velocity},
  journal = {Systems \& Control Letters},
  volume  = {56},
  number  = {3},
  pages   = {245--252},
  year    = {2007}
}

@article{KK:2024,
  author  = {Karafyllis, I. and Krsti{\'c}, M.},
  title   = {Deadzone-adapted disturbance suppression control for global practical {IOS} and zero asymptotic gain to matched uncertainties},
  journal = {Systems \& Control Letters},
  volume  = {185},
  pages   = {105746},
  year    = {2024}
}

@article{PJ:2006,
  author  = {Pepe, P. and Jiang, Z.-P.},
  title   = {A {L}yapunov--{K}rasovskii methodology for {ISS} and {iISS} of time-delay systems},
  journal = {Systems \& Control Letters},
  volume  = {55},
  number  = {12},
  pages   = {1006--1014},
  year    = {2006}
}

@article{JW:2002,
  author  = {Jiang, Z.-P. and Wang, Y.},
  title   = {A converse {L}yapunov theorem for discrete-time systems with disturbances},
  journal = {Systems \& Control Letters},
  volume  = {45},
  number  = {1},
  pages   = {49--58},
  year    = {2002}
}

@article{ROKB:2026,
  author  = {Rodrigues, V. H. P. and Oliveira, T. R. and Krsti{\'c}, M. and Ba{\c s}ar, T.},
  title   = {Distributed event-triggered {N}ash equilibrium seeking},
  journal = {Automatica},
  volume  = {192},
  pages   = {113164},
  year    = {2026}
}

@misc{ROKA:2026,
  author        = {Rodrigues, V. H. P. and Oliveira, T. R. and Krsti{\'c}, M. and Allg{\"o}wer, F.},
  title         = {Event-triggered discrete-time multivariable extremum seeking systems},
  year          = {2026},
  howpublished  = {arXiv preprint arXiv:2607.25223},
  eprint        = {2607.25223},
  archivePrefix = {arXiv}
}

@misc{ROKT:2026,
  author        = {Rodrigues, V. H. P. and Oliveira, T. R. and Krsti{\'c}, M. and Tabuada, P.},
  title         = {Event-triggered {N}ewton extremum seeking for multivariable optimization},
  year          = {2026},
  howpublished  = {arXiv preprint arXiv:2601.14416},
  eprint        = {2601.14416},
  archivePrefix = {arXiv}
}

\end{document}